\documentclass{notices}
\usepackage{xcolor, amsfonts,amssymb,amsmath,amscd,graphicx}

\usepackage{hyperref}

\hypersetup{
  colorlinks=true,
  urlcolor=blue
}

\newcommand{\pullquote}[1]{%
  \par\medskip
  \noindent
   \vskip 0.4\baselineskip
  \noindent\makebox[\columnwidth][c]{%
    \begin{minipage}{0.88\columnwidth}
      \centering
      \small\sffamily\itshape
      ``#1''
    \end{minipage}%
  }
  \par
  \vskip 0.5\baselineskip
  \noindent

  \par\medskip
}

\title{Come for the vibe, stay for the math\\
}

\author{
  Hugo Parlier
  \affil{
    Hugo Parlier is a full professor of mathematics at the University of Fribourg.
    His email address is \texttt{hugo.parlier@unifr.ch}.
  }
  \and
  Bruno Teheux
  \affil{
    Bruno Teheux is an assistant professor of mathematics at the University of Luxembourg.
    His email address is \texttt{bruno.teheux@uni.lu}.
  }
}

\begin{document}

\maketitle


\noindent{\it Our outreach journey began, as many good things do, with a rejection.}\\

Neither of us had started our careers with any notion of becoming math outreachers. But for a decade now, we've shared the creative and exploratory nature of mathematical research far and wide at festivals and expos and schools.  Along the way, we've refined our motivations and developed a practice, of sorts, that we hope others might find inspirational.

In 2017, as newly hired colleagues at the University of Luxembourg, we submitted a proposal for a mathematical exhibit at the Luxembourg Science Festival. The response was polite but unambiguous: no thank-you. Mathematical proposals at this festival were relatively rare, and it wasn't clear to the organizers how one could make our proposed subjects alluring in such a context with lasers, robots and volcanoes. And in fact, the organizers had decided in favor of {\it another} mathematical workshop featuring magic tricks performed by high school students. We expressed our disappointment, and the organizers graciously offered us a free spot when there was a cancellation.

That first festival was an eye-opener for us. Determined to make our exhibit a success, we aimed, with more conviction than experience, to create a space with a lounge type atmosphere, where visitors could enjoy the intuitive and exploratory aspects of mathematics without its formalism, and in contrast to its drier reputation. We presented a collection of games and puzzles, some of our own design, some borrowed, all chosen because they had mathematics quietly in the background. No formulas on the walls, no labels announcing that this was mathematics. Just games, and a willingness to interact.

\begin{figure}[tbp]
\centering
\includegraphics[width=0.8\linewidth]{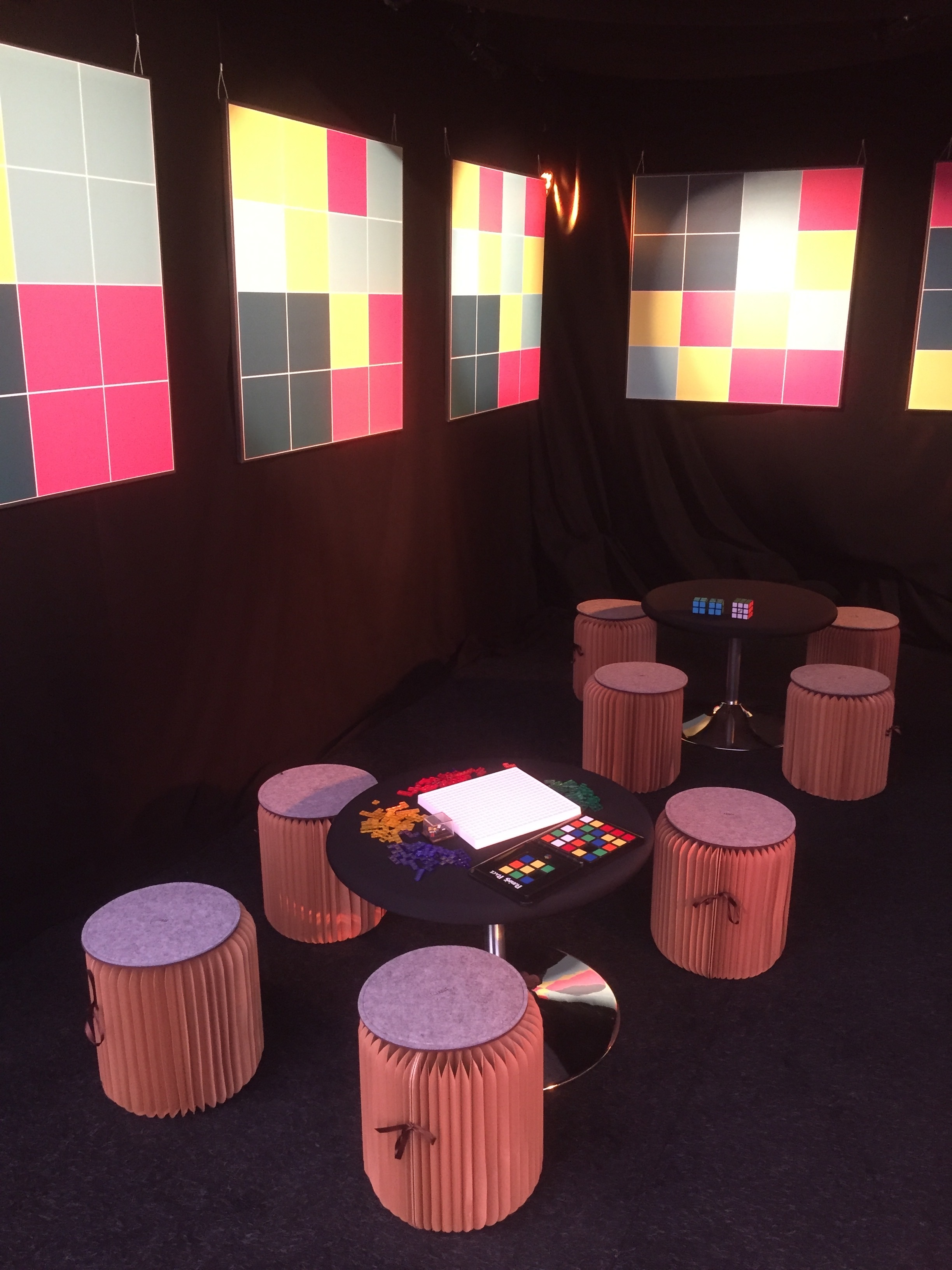}
\caption{Our first setup in 2017}
\end{figure}

What happened that day shaped everything that followed. People sat down, played, argued, shared their joy or surprise, came back to finish a puzzle after having walked away. Children pulled grandparents by the hand to come and try. Adults who had declared themselves irredeemably bad at mathematics found themselves absorbed, twenty minutes later, in solving puzzles. Amid the noise and energy of a festival with afore mentioned lasers, robots, and a working volcano, the lounge was a haven of tranquility. Visitors came for the vibe, and stayed for the fun. 

The exhibit was a success by any measure, and it secured our place among the science outreach community in Luxembourg. More importantly, it gave us an early glimpse of something we have since come to believe more firmly: Encounters between mathematics and wide audiences can produce genuine positive emotions; and these emotions are what make such encounters memorable. While some understanding of why it worked was immediate, it has deepened and evolved with each subsequent event.


Emboldened, we began creating exhibits systematically. We proposed activities for science festivals in subsequent years, each time with new themes and material. The titles give a sense of the range: \textit{Jouer, penser, r\'esoudre (Play, think, solve)}, \textit{Unpuzzling Mathematics}, \textit{Mathematics in Black and White}, \textit{The Simplicity of Complexity} , \textit{Reshape},  \textit{Traces}. For each of these, we created and designed activities from scratch: games, puzzles, and visual explorations rooted in genuine research questions, but accessible through simple, immediate rules.

\section*{And where is the math?}

Here is a sampling of the type of activities we present at events\footnote{See \url{http://mathword.lunot.eu} for a short account of these activities and a recent list of events we have contributed to}.

Outreach activities often have a combinatorial flavor, and many of ours are no exception. We also try to bring other mathematical aspects into play, many of them related to our own research, which include geometry, group theory and topology. 

\begin{figure}[tbp]
\centering
\includegraphics[width=0.8\linewidth]{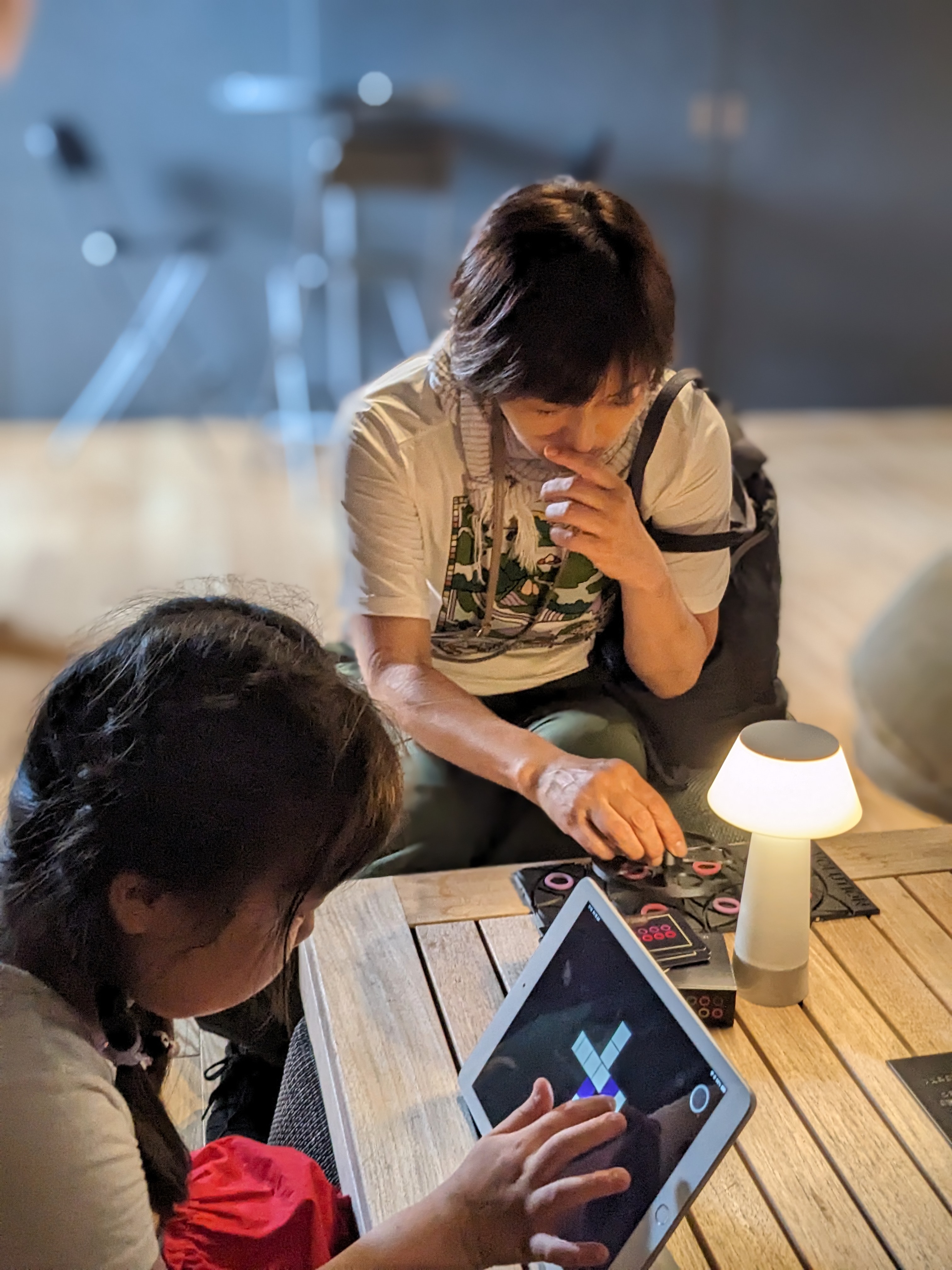}
\caption{Bridging generations through games}\label{fig:bridge}
\end{figure}

\begin{figure}[tbp]
\centering
\includegraphics[width=0.8\linewidth]{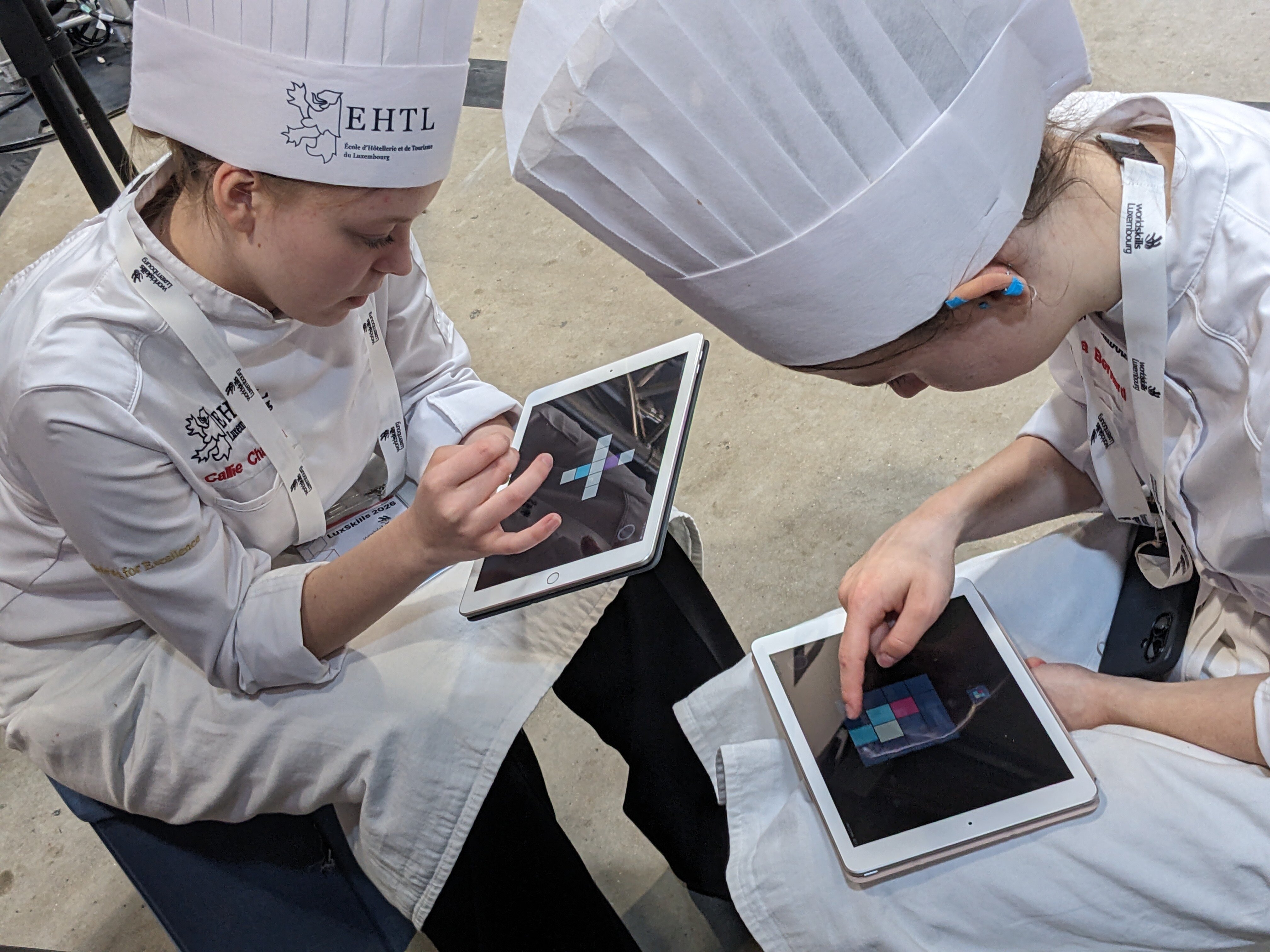}
\includegraphics[width=0.8\linewidth]{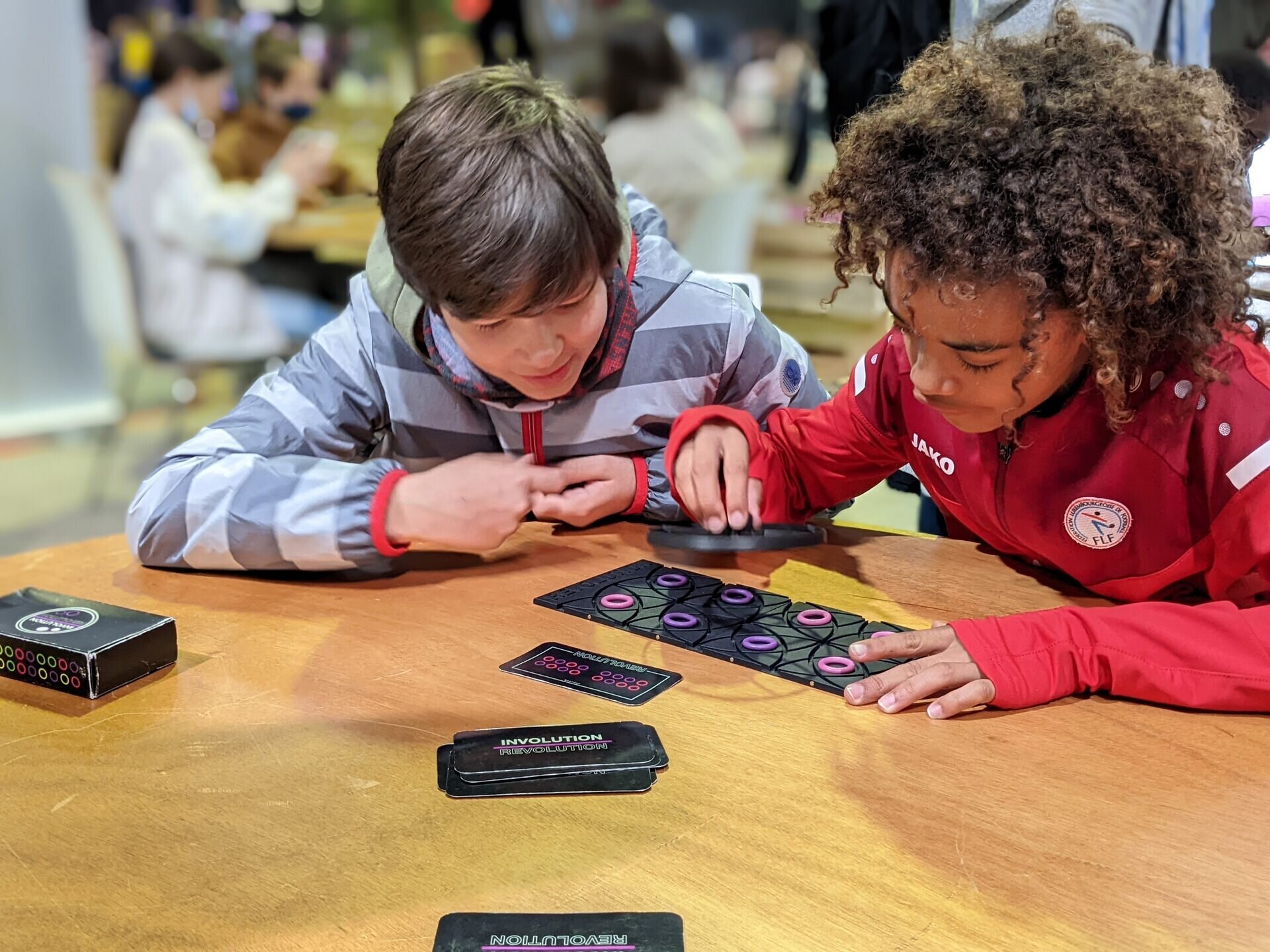}
\caption{Quadratis and Revolution in action}\label{fig:quadrevo}
\end{figure}

\textit{Quadratis}\footnote{Quadratis was co-created by H. Parlier and P. Turner with development by M. Gutierrez and R. Juarez and is available on http://quadratis.app.} is a puzzle type game in which the player is presented with a pattern of colored square tiles, which are then shuffled and your goal is to restore the original pattern by sliding tiles (see Fig.~\ref{fig:bridge} and \ref{fig:quadrevo}).

The shuffling and sliding are subject to rules determined by the way the squares are glued together thus, in effect, by the topology of an underlying surface. Sides may be identified in unexpected ways, so that a tile sliding off the right edge of the board may reappear at the bottom left. What looks like a familiar sliding puzzle turns out to encode the geometry of so-called square translation surfaces, flat surfaces built by identifying edges of polygons and objects of active research in geometry and dynamical systems. In effect, by playing the game, you are exploring types of moduli spaces, vast structures which encode different geometries that surfaces can take. 

The puzzles range from straightforward to genuinely difficult. They invite deeper questions: How many moves are needed to reach a given configuration? Is any configuration reachable from any other? What does an optimal algorithm look like? And these questions are often open research problems, some closely related to mysteries about the underlying moduli spaces. The game is freely available for download, which means that the encounter can last long after the expo.

\textit{Revolution} and its companion \textit{Involution}\footnote{Revolution and Involution were created and designed by H. Parlier and B. Teheux.} are physical reconfiguration puzzles, distant cousins of the Rubik's cube, though more approachable. Difficulty increases gradually rather than all at once (see Fig.~\ref{fig:quadrevo}). Moves consist of rotating a wheel to rearrange a sequence of colored rings and the goal is to achieve a prescribed pattern. The games were again created from scratch from the starting idea of being as simple as possible. The fact that you rotate 4 rings at a time is linked to us seeing what we could come up with by just playing with our fingers. The games have immediate links to subgroups of permutation groups, but are in fact inspired by more topological aspects, such as braids.

Part of the ``showcasing" the mathematics behind the activities involves making the state space of the puzzle visible: we display a graph whose nodes are configurations and whose edges are legal moves, so that visitors can see, literally, the structure of the problem they are navigating. This graph aids in the quest to bring visitors to the frontier of knowledge. The natural questions --- can I always get from here to there, and how efficiently? --- are again open problems in graph theory and combinatorics. The challenges themselves have developed over time, based on tests, experience, and cognitive aspects, but also with our understanding of the underlying graphs. We found that using more sophisticated invariants such as treewidth is useful for understanding, and hence creating, more intricate challenges.

\begin{figure}[tbp]
\centering
\includegraphics[width=0.95\linewidth]{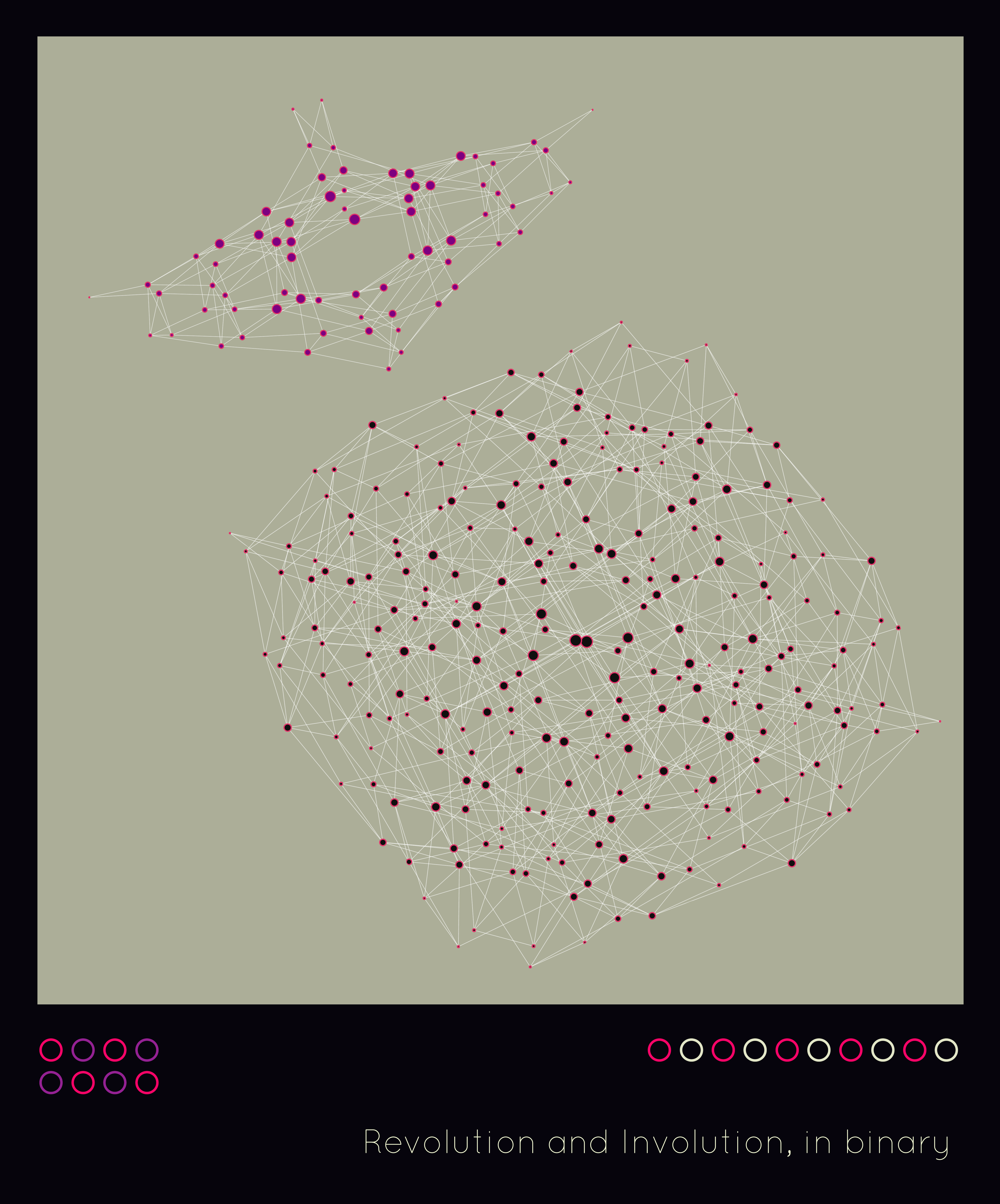}
\caption{State space of the Revolution and Involution games}
\end{figure}

One of the recurring surprises for visitors, and one of the things we most enjoy pointing out, is that we cannot always solve the puzzles we have designed. The rules are ours; the questions the rules raise are not. This is, of course, precisely the situation of the working mathematician: you define the object, and then discover that you do not understand it. Visitors find this genuinely astonishing. Some of the best moments are when a visitor, having wrestled with a puzzle for some time, asked a question to which the only honest answer was: ``We do not know. And in fact we suspect that --- at this point in time --- no-one knows."

\pullquote{The rules are ours; the questions the rules raise are not.}

\textit{Life Lines}\footnote{Life Lines, created by H. Parlier and B. Teheux, was tried out in early versions during a second residency at Expo 2020 Dubai in 2022, and first showcased in its current form at the \href{https://epfl-pavilions.ch/en/exhibitions/shapes}{Shapes} exhibit at EPFL Pavilions in January 2025. See also \href{https://lifelines.cloud/}{https://lifelines.cloud/}.} is of a different character. Visitors are invited to draw a curve with their finger on a touchscreen, subject only to the constraint that it travels from left to right. The drawing is then colored automatically, in a way that seems at first mysterious. Different regions of the drawing receive different colors, according to a rule that is precise but not immediately apparent. The coloring is determined by the winding number of the curve with respect to each region it delimits, a concept from topology that records how many times, and in which direction, the curve wraps around each enclosed area.

\begin{figure}[tbp]
\centering
\includegraphics[width=0.8\linewidth]{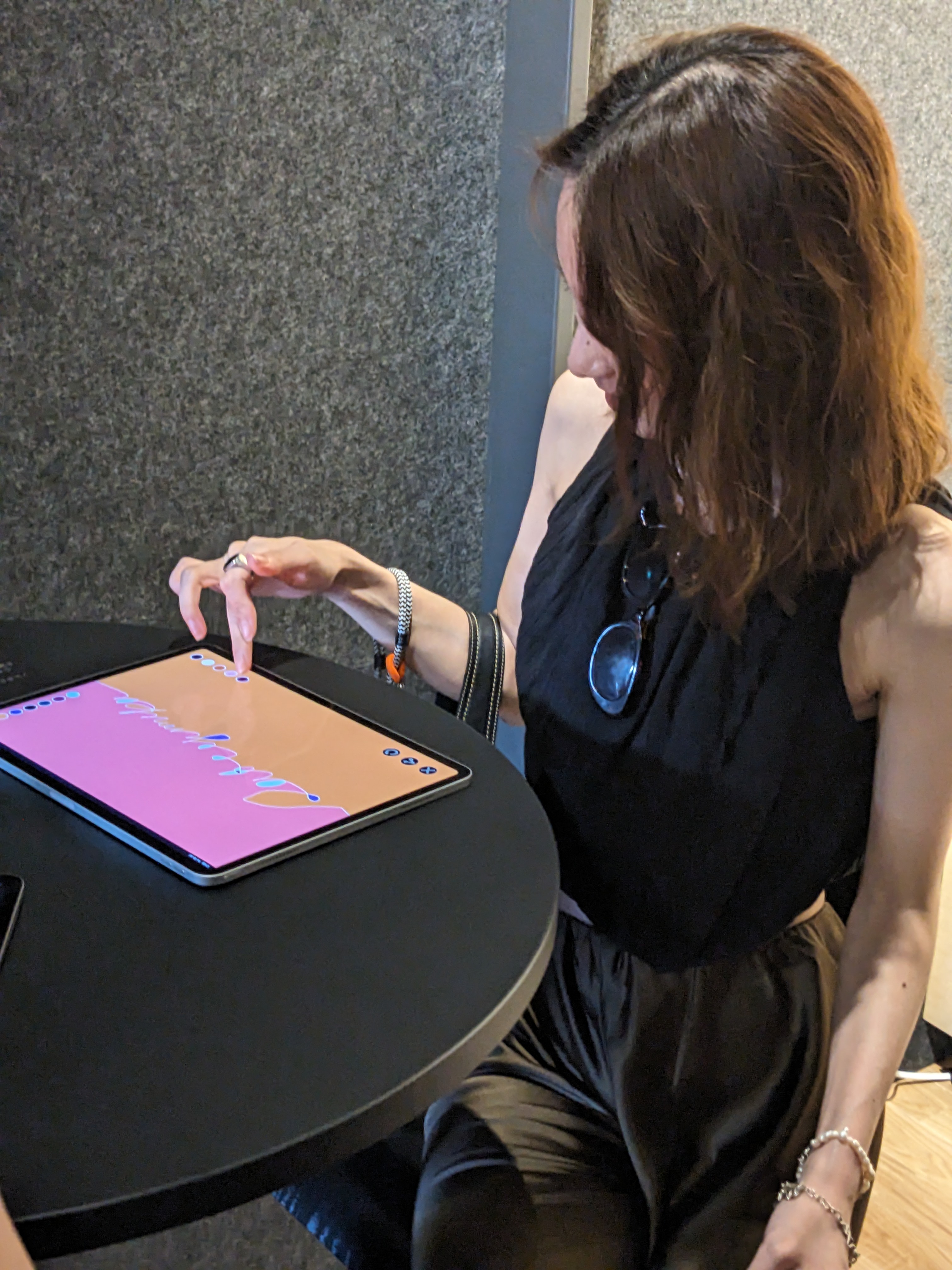}
\caption{Drawing with Life Lines}
\end{figure}

The colored drawings are displayed in a continuously evolving shared installation. Visitors can browse the drawings collected, compare their creations to those of other visitors, and watch morphings --- continuous deformations of one drawing into another, with the colors shifting in real time. There is structure underlying these deformations, a geometry on the space of drawings itself. This is similar to some of the previously described activities; each drawing is a point in an underlying, infinite dimensional, moduli space of all possible drawings. To measure distance between drawings, we use the so-called Fr\'echet distance, informally known as the dog-leash distance. It measures proximity between curves as paths rather than as point sets, and is related to active research in computational geometry. We also use it to compute averages of curves, which allows us to explore, in a precise way, what the aggregate of all visitor drawings looks like.

\begin{figure}[tbp]
\centering
\includegraphics[width=0.8\linewidth]{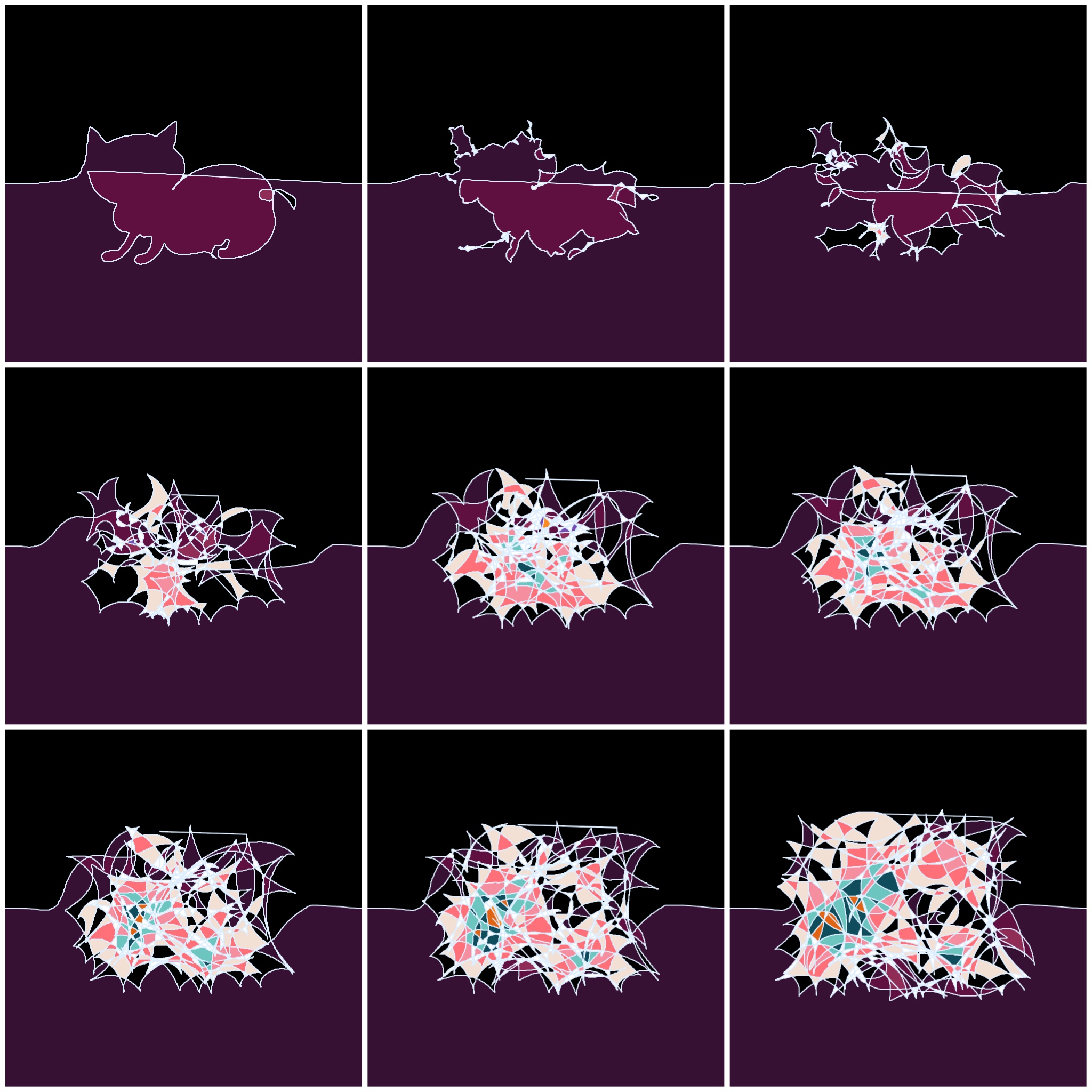}
\caption{A morphing between figurative and abstract drawings}
\end{figure}

Importantly, \textit{Life Lines} offers something different from the puzzle activities: an invitation to create rather than to solve. Some visitors draw abstract patterns while others draw recognizable figures --- animals, faces, portraits, imaginary landscapes. All of them were surprised by the colors that appeared. Many came back to draw again, trying to understand the logic, trying to produce a particular effect, trying to make something beautiful. The contrast with the puzzle activities is deliberate: some visitors want to solve, others want to create; some want to compete against a puzzle, others want to collaborate; some want to work alone while others come alive in groups. By offering multiple modes of engagement, we are able to speak to an even wider range of people and, incidentally, to showcase some of the diversity of mathematics itself.

\pullquote{By offering multiple modes of engagement, we are able to speak to an even wider range of people and, incidentally, to showcase some of the diversity of mathematics itself.}

Over the course of various presentations --- in Dubai, Lausanne, Osaka --- tens of thousands of drawings were collected; a single Dubai residency alone yielded 20,000 drawings and 50,000 accompanying puzzle solutions. These numbers meant something beyond mere statistics: they represented an extraordinary range of people who had, for a few minutes or a few hours, engaged seriously and joyfully with mathematical ideas. The collection is not primarily remarkable as a geometric object, though the mathematics it involves is genuine and current. It is remarkable as an archive of human creativity shaped by a mathematical rule: a record of what people choose to draw when given a blank screen, an automatic coloring, and no further instruction.

\begin{figure}[!htbp]
\centering
\includegraphics[width=1\linewidth]{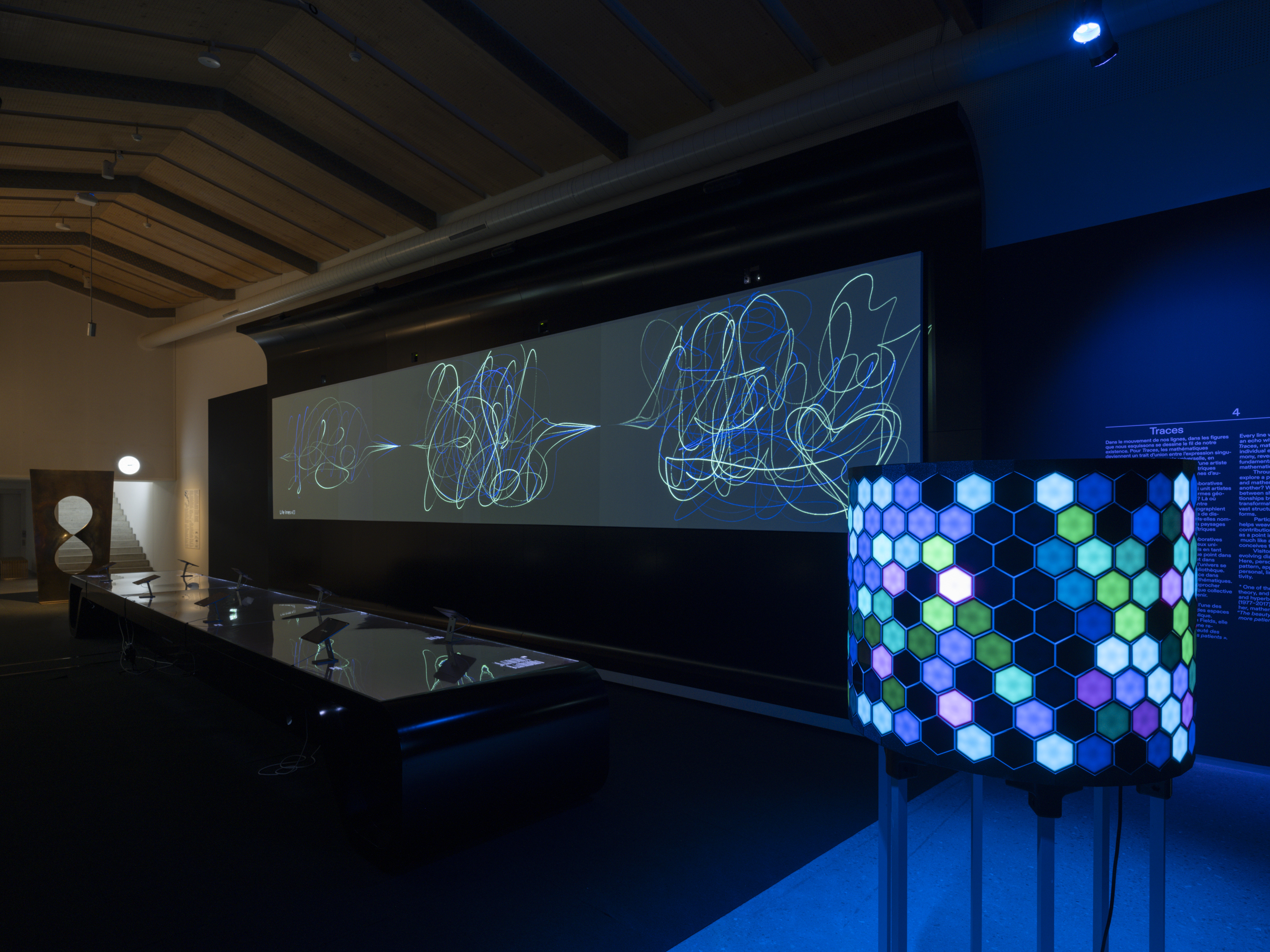}
\caption{The Shapes exhibition at EPFL Pavilions, 2025}
\end{figure}

\begin{figure*}[t]

\fbox{%
  \begin{minipage}{0.98\textwidth}
  \textbf{The cyclicity cylinder.}

  \medskip
  One of the objects we built has an unusual story, even by our standards.

  \medskip

  This 90cm diameter cylinder stands at roughly 2m, and the upper part is tiled by touch
  screen hexagons whose color and light intensity is dictated by a CPU in the
  background. We conceived the object itself {\it before} having a clear idea of
  what we wanted to do with it.

\begin{center}
\includegraphics[height=6cm]{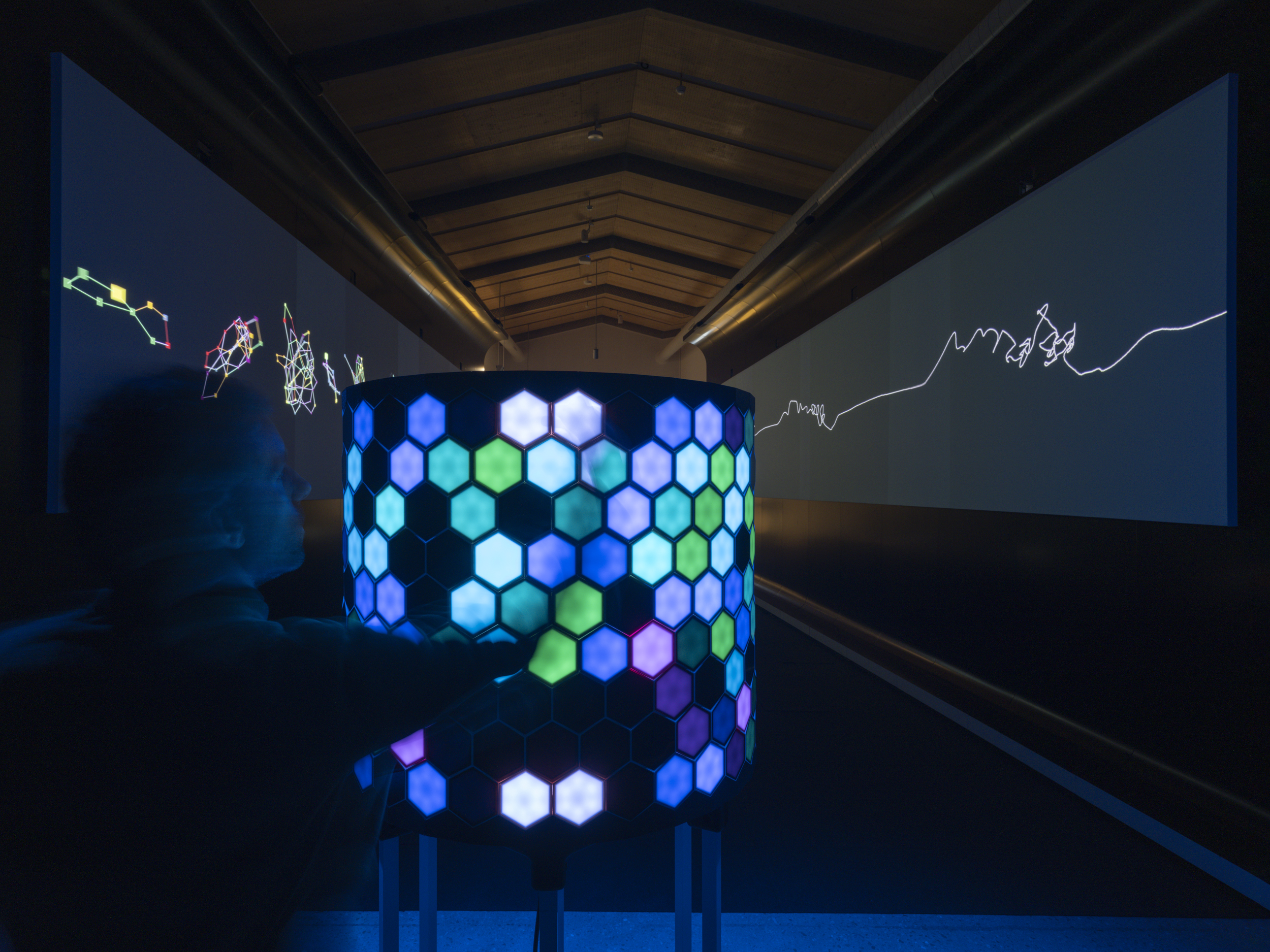}
\hfill
\includegraphics[height=6cm]{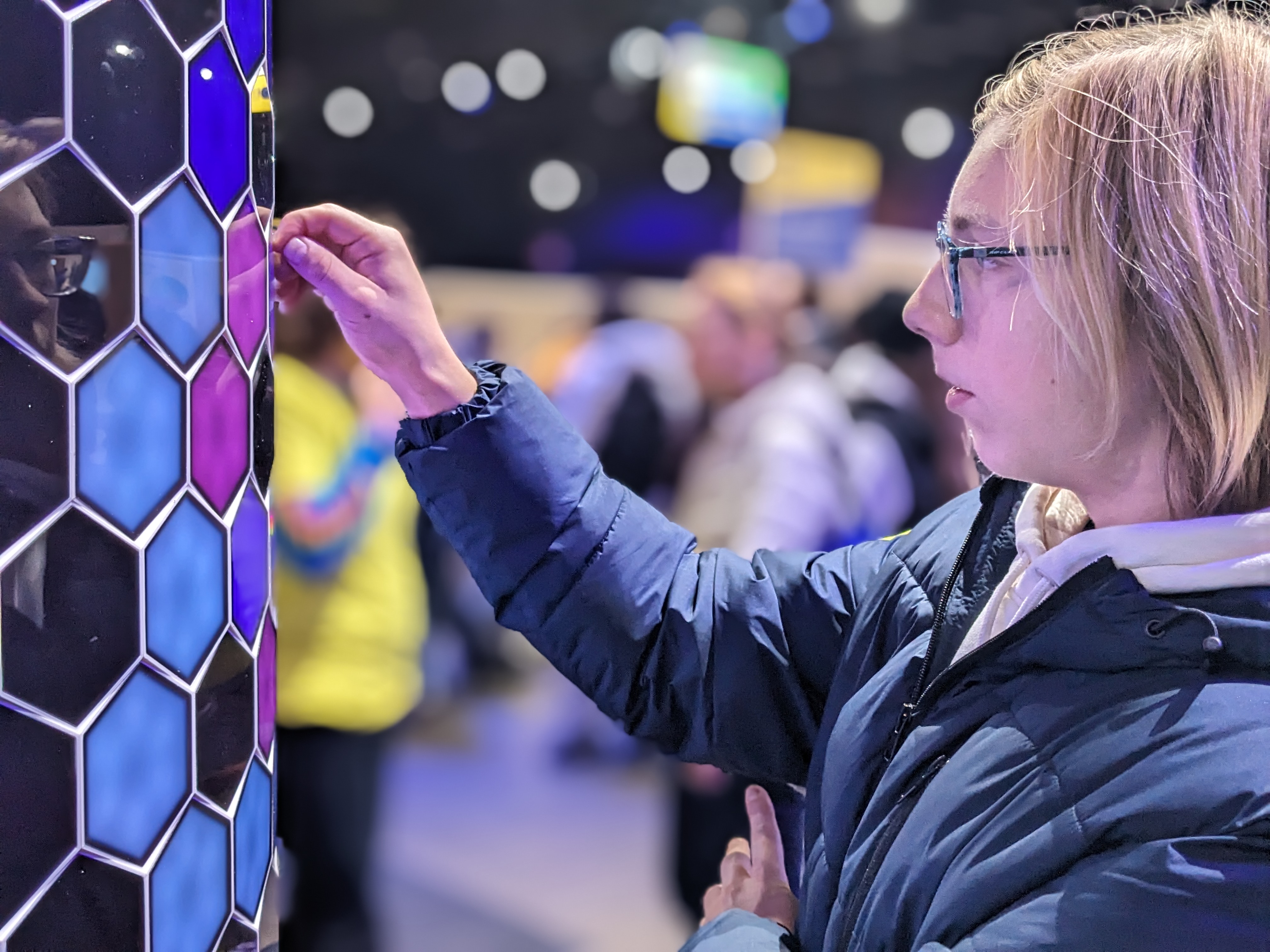}
\end{center}

With the cylinder as inspiration, we created the game {\it Cyclicity}, ultimately based on homology, which can be enjoyed on the cylinder as a playing board or on touch screen tablets. The game was featured on {\it Take Off}, a Luxembourgish science TV show (here is a \href{https://www.youtube.com/watch?v=UO1Buyh_Fmk&t=1470s}{link} to the episode). We also turned the cylinder into an interactive sculpture featured at
the {\it Shapes} exhibition, which ultimately represents a four dimensional
torus. In both cases, we only recognized the mathematics behind the activities
in retrospect. And in complete honesty, we are not completely sure how, when or
why we came up with the idea of building the cylinder in the first place.  \end{minipage}%
}
\end{figure*}

\section*{Learning by doing}

The decision to make everything ourselves was not, at first, entirely deliberate. Over time, it developed into a firm conviction. When you create your own material, it can be precisely tailored to the experience at hand. You decide what the entry point is, how the difficulty progresses, and what the visitor is invited to wonder about. The design is part of the  atmosphere, from the texture of objects, their colors, to the way a puzzle sits on a table and invites a hand. We are not designers by training, but we found that caring about these details made a difference.

There is another, less obvious, advantage. When an experience is designed before the mathematics are fully worked out, the activity naturally generates genuine open problems. The rules of a game may be simple enough for anyone to begin playing immediately, the questions they raise, about reachability, optimality, structure, may be entirely open.
This is not a shortcoming but a feature. It means the activity sits at the boundary of knowledge, and that the mediator\footnote{The words "mediation" and ``mediator" are less used in this context in English, in favor of "outreach" and ``outreacher" but they illustrate a vital part of the process which involves mediating between the topic and suspicious visitors.} can honestly say, \textit{I don't know!} in response to a visitor's question. It also means that the same activity speaks to visitors at very different levels. A novice will simply jump in and experience the pleasure of attempting a puzzle while a trained mathematician might see the configuration space and wonder about its structure. Over the years, our designed outreach activities have generated multiple mathematical questions that interest us as researchers.

Crowds at festivals or large-scale events provide a demanding but invaluable testing environment. With hundreds of visitors passing through each day, you quickly learn what works and, more importantly, what does not. You start to feel the moment when a visitor's attention wanders and begin to understand why. You learn to distinguish the visitor who wants a hint from the one who needs to be left alone to think. We always introduce new material alongside proven activities, mitigating the risk of failure and pinpointing benefits and shortcomings of the new material.

We also learned something more humbling: that even the most carefully designed activities require improvisation and revision. As much work as we put into creating instructions, the first time a new activity meets real visitors, something unexpected always happens. A step that seemed transparent turns out to confuse. A puzzle we judged to be of medium difficulty proves to be either too easy or completely impenetrable. The first run of anything new involves rapid, real-time negotiation between what was planned and what is actually happening in front of you. It is stressful. It is also consistently entertaining. The situations that arise from these mismatches are the ones we recount to each other years later.

\section*{Expanding our comfort zone}

Science festivals draw a naturally curious audience. The visitors arrive, in some sense, prepared to be intrigued. A game festival is another matter entirely. When we were invited to participate in the \textit{Game On} festival in Luxembourg City, we found ourselves competing not with volcanoes but with wonderful new board games with beautiful production values, devoted communities, and years of refinement behind them.

What we came to think of as our ten-second rule proved its worth here. Many excellent board games require fifteen minutes of rulebook reading before they can begin. Our activities required a gesture: watch once and then try. Visitors who might never have paused at a
mathematics exhibit are hooked by the simplicity of the entry, and surprised by the depth they found inside. The audience itself was different: families and gaming enthusiasts with no particular attachment to mathematics or to science, who had come simply to play.

\pullquote{Our motto, which was only half a joke, became `Come for the slide, stay
for the math'.}

The same period saw us brought into contact with an entirely different kind of venue. Through a series of circumstances that we could not have anticipated, we were invited to present some of our activities at the Luxembourg Pavilion at Expo 2020 Dubai. The pavilion itself had been designed in the shape of a M\"obius band, which we took to mean a sign of destiny, and its main attraction was a giant slide that turned out to be a visitor magnet. Our motto, which was only half a joke, became ``come for the slide, stay for the math".

It worked. With something on the order of a thousand active visitors a day, the scale of engagement was unlike anything we had experienced before. The audience was genuinely international: visitors from across the world, of all ages, arriving with every level of prior engagement with mathematics ranging from deep enthusiasm to profound indifference. Most left with a smile, perhaps slightly puzzled, and certainly with at least a slightly altered view on the subject. The pavilion asked us to return for a longer second residency, during which we also tried out early versions of what would become \textit{Life Lines} (described above).

Expo 2025 Osaka presented challenges of a different order: a two-week residency at the Luxembourg Pavilion, plus events we organized directly with the Swiss and Belgian Pavilions --- three venues in a country whose language and cultural context were quite different from anything we had previously encountered.

We were warned, by people who knew the context well, that engaging Japanese visitors without fluent Japanese-speaking mediators would be very difficult. We arranged for Japanese-speaking colleagues (not necessarily mathematicians) to attend and be available to assist. We adapted our visual materials. We thought carefully about which activities could be demonstrated most clearly without words.

The rules for most of our activities are sufficiently simple and sufficiently visual to be conveyed by demonstration, via a gesture, a move, an invitation to try.  The mathematics does not require verbal explanation to be entered. And once a visitor begins, the experience speaks for itself. By sharing something new, surprising and which brings joy, allows to create an emotional bond that transcends cultural and language barriers.

One aspect of the Osaka audience that we had not anticipated, and that moved us considerably, was the depth of commitment brought by many visitors to the activities. Day after day, visitors returned and spent hours working through every puzzle we had laid out. Some as young as six, others well beyond retirement age. A six-year-old solving a puzzle that had defeated every adult nearby produces a reaction that needs no translation. The delight is immediate, shared, and completely recognizable. This commitment sometimes took concrete form: On several days, visitors queued for ninety minutes in order to return to activities with which they had already spent time.

We were also helped, no doubt, by the particular Japanese appreciation for games and for the pleasures of careful thought. But we would be cautious about attributing too much to cultural specificity. The depth of engagement we saw in Osaka was not categorically different from what we had seen in Europe or in Dubai, it was just a more concentrated version of the same phenomenon.

\section*{What we learned}

Across festivals, expos, school visits and with activities now reaching tens of thousands of visitors, certain convictions have settled into place. The diversity of the audience is not an incidental feature of what we do, it is one of the main reasons why we do it. We have worked with children encountering puzzles for the first time and with elderly visitors bringing decades of games experience. With families who knew nothing about mathematics research and with professional mathematicians who wanted to discuss the open problems embedded in the activities. With visitors with whom we shared no languages and visitors who came to converse. With participants full of curiosity, and others with indifference, but who always left with an altered point of view. 

Such reach has been possible, we believe, because the goal has never been to transmit mathematical knowledge. It has been to create a positive emotion about mathematics, one that is genuine, immediate, and, we hope, lasting. An emotion that precedes understanding and, if it takes hold, may eventually motivate understanding. Pleasure is not a wrapper around the mathematics --- it is the form in which the mathematics is first encountered and experienced. This approach also reaches people who are entirely unmoved by the argument that mathematics is useful, an argument that, however well-founded, does not resonate with everyone.

This conviction has practical consequences that are sometimes counterintuitive. It means that the mathematical depth of an activity is not in tension with its accessibility. On the contrary, the deeper the mathematics, the richer the activity tends to be at every level. It means that the atmosphere and the design are not decorative but substantial --- they are what determines whether an encounter is possible at all. And it means that what the visitor takes away is not, primarily, information about mathematics, but an experience of mathematical thinking, which is something else entirely.

We have also learned how important it is to be honest about what outreach cannot do. A visitor who spends twenty minutes with one of our activities does not leave with a new mathematical competency. What they leave with, if the encounter goes well, is a memory, increased curiosity, and perhaps a slightly revised picture of what mathematics is and what mathematicians do. That is not nothing. One could argue that it is, in fact, quite a lot. But it requires the mediator to keep to their modest objectives, and to resist the temptation, familiar to anyone with a teaching instinct, to explain more than the moment requires.

One small episode illustrates the depth of what can happen more concretely than any statistics. When the time came to curate the drawings collected during our Dubai residencies for a data-sonification project \textit{The Sound of Data}\footnote{The Sound of Data project was a joint project of the Luxembourg National Science Fund, the Rockhal, the LIST and the University of Luxembourg and was part of Esch 2022, when Esch sur Alzette was the European capital of culture.}, we sat down with the archive for the first time with a view to selection rather than collection. What surprised us was that for a significant number of the drawings, we could remember the person who had drawn them. The process brought back faces, memories of conversations, and moments of surprise or delight. We had not realized, while they were happening, how thoroughly these encounters had affected us. The archive that looked like a dataset turned out to be something closer to a collection of memories.

These moments and memories are deeply fulfilling. We mean this in a specific sense. Moments of joy for visitors, especially for a non-captive audience who elected to come, are of course moments of joy for us. There is a particular quality to the moment when a visitor who was about to leave becomes, instead, absorbed. Something shifts: the body language, the attention, the relationship to the activity. There is also the visitor who, having finally solved a mystifying puzzle, comes back looking for a high five, with a brightness in their eyes that is unmistakable. We have watched moments like these happen thousands of times, across venues and with a great diversity of visitors, and yet they never become routine. If anything, the accumulation of these
moments has made us more attentive to what they mean. Something genuine is being transmitted. And it is not a theorem, not necessarily a concept, but something about what it feels like to encounter and resolve a genuinely hard problem. This transmission is felt on both sides.


 \section*{Principles we (sometimes) follow}

The principles below were not written down before any of the activities described above took place. They emerged, gradually, from the experience of standing in front of thousands of visitors and observing what created connections and what did not. There are an arbitrary number (10), mostly to ensure that we can still count them (on our fingers). 

They are design patterns rather than prescriptions, tendencies we now follow largely by instinct. In particular, one thing they are {\it not} is rules. On occasion we overlook these principles, sometimes on purpose.  They are more of a record of what experience has taught us, in the hope that they may be useful to others starting down a similar path.\\



\noindent{\bf 1. Multiply entry points}\\

No one should be written off in advance, and the first challenge is simply getting someone to stop and try. People take that first step for very different reasons.
\pullquote{...the first challenge is to get someone to stop and try.}
Maybe they like games, maybe the colors caught their eye, or maybe a friend is already playing. Maybe the atmosphere feels welcoming because the mediator looks approachable, or because the objects on the table are beautiful. Maybe someone just wanted to sit down and then felt obligated to at least try the activity. None of these reasons are mathematical, and that is precisely the point. But the more reasons there are, the wider the audience. Mathematics comes second, the invitation comes first.\\


\noindent{\bf 2. Keep it simple}\\

The entry point of any activity must be immediate. A visitor at a science festival, a game fair, or a world expo has dozens of competing distractions. You have at most ten seconds to catch their attention. But once caught, attention can last, sometimes far longer than visitors themselves would expect.

\pullquote{Mathematics is not a spectator sport: to enjoy it, you have to engage, and engagement begins the moment the rules are understood.}

We think of this as the ten-second rule: the instructions of any activity must be fully deliverable in under ten seconds, so that the threshold to trying is as low as possible. Mathematics is not a spectator sport: to enjoy it, you have to engage, and engagement begins the moment the rules are understood.\\


\noindent{\bf 3. Be progressive}\\

At an open event, progression in difficulty is a key factor in making the activity inclusive. The challenge increases gradually, with natural stopping points where a visitor can feel satisfied and leave, or feel curious and stay. The early stages must be easy enough to guarantee a first small victory quickly. Once a visitor has genuinely solved something, their willingness to attempt harder things grows considerably. In the best cases, the pleasure of progression becomes something close to a positive addiction.\\

\noindent{\bf 4. Put the fun first}\\

Let pleasure drive the design, yours as much as the audience's. Build an experience that is genuinely enjoyable on its own terms, and let the mathematics appear naturally. It should arrive as a discovery rather than a lesson. The visitor is not receiving mathematics, they are gradually but actively noticing it. \\


\noindent{\bf 5. Be modest}\\

Resist the educator's instinct to explain more than the moment requires. It is already fantastic if visitors leave with a warm memory, a genuine curiosity, a slightly revised picture of what mathematics is. The goal is to create a need to continue exploring, not to satisfy it on the spot. An outreach activity is a teaser. And a good teaser does not give everything away.\\


\noindent{\bf 6. Be holistic}\\

Appeal to all senses. Every aspect of the environment contributes to whether the encounter happens at all. Furniture, lighting, the texture of objects, sound, the visual design of materials, the way an activity is laid out on a table and invites you to jump in. None of this is trivial. Being holistic and paying attention to all these aspects is also a factor of inclusivity. When every sensory dimension of the experience has been considered, the activity becomes accessible to a wider range of people, including those whose relationship to mathematics is primarily visual or tactile rather than symbolic. Design is a tool.\\

\noindent{\bf 7. Go where you are unexpected}\\

Science outreach tends to circulate within its own networks including science festivals, university open days and school visits. These are important and valuable venues, but they largely reach people who were already curious. Venturing outside these circuits means reaching a genuinely different audience, and one that deserves the encounter just as much. At a game festival or a cultural event, mathematics arrives as a surprise. Visitors have no prior expectations, and so cannot be disappointed. In our experience, the impact of showing up where you are not expected is well worth the effort of getting there.\\



\noindent{\bf 8. Love your audience}\\

Don't be afraid to have genuine affection for the people you share your activities with. If you do, you will naturally pay attention to who they are, and that attention will shape everything from how you engage, when you step in and when you step back. Not every encounter produces a spark, but some produce something that stays with both parties long after the event. The sparkles in their eyes will mirror yours. If these moments matter to you, that is your best reason to keep doing this.\\
 

\noindent{\bf 9. Mediation is subtle}\\

Good mediation is like framing a photograph: It is as much about what you leave out as about what you put in. The impulse to explain, to teach, to fill silence  with information, is familiar to anyone trained as an educator. In this context, it is precisely what must be resisted.  A visitor who is struggling needs encouragement, not a solution. Contact matters more than guidance. The mediator's role is to be attentive, and supportive of the visitor's own exploration, not to direct it. Mediation requires training, and this mostly comes through experience. When it works, it creates something stronger than understanding: an emotional bond formed around a shared moment of mathematical experience.\\



\noindent{\bf 10. Seek institutional support}\\

Our outreach presentations required sustained institutional support to exist. For researchers who wish to invest seriously in outreach, alignment with their institution's strategy is essential. This is not primarily for the institution's benefit, but for their own. Close collaboration with communication and logistics teams is equally important. They bring complementary expertise, and they become partners, not service providers. Significant funding is sometimes also required, namely for design, fabrication and development, and for travel. 

In our case, the National Research Fund in Luxembourg has been an essential partner. Without public funding of this kind, our projects would not have been possible.

Our universities, and both the rectorate and the mathematics departments in Luxembourg and Switzerland, have been systematically generous and supportive, including making outreach a part of our professional duties. While all departments might not have the will or the means to dedicate resources to outreach, it is often a good ``deal", as it helps anchor the department in the local community and brings a good deal of visibility. In our case at least, a great amount of personal time and energy go into these projects for the benefit of all involved, including the institution.

 \section*{Final thoughts}
  
While what is explicitly described here is primarily large scale outreach activities, the philosophy behind the approach permeates how we communicate in other contexts. This ranges from how we share what we do with loved ones, how we teach and how we present and even write about research.

We hope our description doesn't make starting in outreach seem daunting. The first steps in outreach rarely begin with large scale activities. Start small and local and don't hesitate to piggy back on existing events. Try ideas on friends and family or anyone you can get ahold of who will give you genuine feedback. Finding out that an activity doesn't work in front of an elementary school class can be terrifying. There are many ways to share math, but the ultimate test is in front of an audience that is not already sold on the activity.  Working with collaborators, namely those that you can shape and share a common vision with, is rewarding and fulfilling. You don't have to begin by creating activities from scratch: repurpose or adapt anything that fits your vision. And just like in math, if a promising idea doesn't work, don't hesitate to drop it. 

Just as mathematics need not be justified only by its usefulness, outreach need not be justified only by its importance. A common argument for encouraging outreach is that it matters (because it does!).  A perhaps equally compelling argument is that doing outreach is a lot of fun. It can be a welcome change from the sometimes solitary pursuit of research, and the sharing of ideas with only a handful of specialists. It is fulfilling, and it might even change the way you think and communicate about math. 

 This article is the fruit of our very personal journeys, and anything that vaguely resembles a rule on ``how one should do outreach" is just a rule that we use for personal guidance. As everything, outreach benefits immensely from variety and diversity, and this includes with respect to underlying principles or philosophies. We very much hope there are fellow outreachers and/or mathematicians who disagree strongly with us as it is a sign of variety.

\medskip

\noindent{\bf Acknowledgements.}

Many of the ideas of this article were first presented in an article for the Swiss VSH-AEU-Bulletin entitled ``The simplicity of complexity, a story of mathematical outreach", and at different talks, including an invited talk at the Joint Math Meetings of the AMS in January 2026. We thank the organizers for the opportunity to share our experience.

\end{document}